\documentclass[12pt]{article}
\usepackage{graphicx}
\usepackage{amssymb}
\usepackage{amsmath}
\numberwithin{equation}{section}

\begin{document}
\author{Lev Sakhnovich}
\date{March 01 , 2007}
\textbf{Explicit Rational Solutions of Knizhnik-Zamolodchikov
Equation.}
\begin{center} Lev Sakhnovich  \end{center}
735 Crawford ave., Brooklyn, 11223, New York, USA.\\
 E-mail address: lev.sakhnovich@verizon.net
\begin{center}Abstract \end{center}
We consider the Knizhnik-Zamolodchikov system of linear differential
equations. The coefficients of this system are rational functions
generated by elements of the symmetric group $S_{n}$. We assume that
parameter $\rho=\pm{1}.$ In previous paper [5] we proved that the
fundamental solution of the  corresponding KZ-equation is
rational. Now we construct this solution in the
explicit form.
\newpage
\section{Introduction }
1.We consider the Knizhnik-Zamolodchikov differential system (see
[3])
\begin{equation}
\frac{dW}{dz}={\rho}A(z)W,\quad z{\in}C,\end{equation} where $A(z)$
and $W(z)$ are $n{\times}n$ matrices. We suppose that $A(z)$ has the
form
\begin{equation}
A(z)=\sum_{k=1}^{n-1}\frac{P_{k}}{z-z_{k}}, \end{equation} where
$z_{k}{\ne}z_{\ell}$ if $k{\ne}\ell$. The matrices $P_{k}$ are
connected with matrix representation of the
symmetric group $S_{n}$ and are defined by formulas (2.1)-(2.4).\\
\textbf{Remark 1.1.} Equation (1.1) is one of the equations
belonging to the consistent KZ-system. In section 3 we shall consider
 the whole KZ-system.\\
 A.Chervov and D.
Talalaev formulated the following interesting conjecture [2].\\
\textbf{Conjecture 1.1.} \emph{The Knizhnik-Zamolodchikov system
$(1.1),(1.2)$ has a rational fundamental matrix solution when
parameter $\rho$ is integer.}\\
We have proved this conjecture [5] for the case when $\rho={\pm}1$
In the present paper we solve
the corresponding KZ-equation in the explicit rational form.
 In the case $S_{4}$ the explicit solution was constructed in
 paper [6].\\
 2.  In a
neighborhood of $z_{k}$ the matrix function $A(z)$
 can be represented in the form
\begin{equation}
A(z)=\frac{a_{-1}}{z-z_{k}}+a_{0}+a_{1}(z-z_{k})+...,\end{equation}where
$a_{k}$ are $n{\times}n$ matrices. We investigate the case when
$z_{k}$ is either a regular point of $W(z)$ or a pole. Hence the
following relation
\begin{equation}
W(z)=\sum_{p{\geq}m}b_{p}(z-z_{k})^{p},\quad b_{m}{\ne}0
\end{equation} is true. Here $b_{p}$ are $n{\times}n$ matrices.
We note that $m$ can be negative. \\ \textbf{Proposition
1.1.}(necessary condition, (see [4])  \emph{If the solution of
system $(1.1)$ has form $(1.4)$ then $m$ is an eigenvalue of
${\rho}a_{-1}$.}\\
We denote by M the greatest integer eigenvalue of the matrix
${\rho}a_{-1}$. Using relations (1.3) and (1.4) we obtain the assertion.\\
\textbf{Proposition 1.2.}(necessary and sufficient condition, (see
[4]) \emph{If the matrix system
\begin{equation}
[(q+1)I_{n}-{\rho}a_{-1}]b_{q+1}=\sum_{j+\ell=q}{\rho}a_{j}b_{\ell},
\end{equation} where $m{\leq}q+1{\leq}M$, has a solution $b_{m},
b_{m+1},...,b_{M}$ and $b_{m}{\ne}0$ then system $(1.1)$ has a
solution
of form $(1.4)$.}\\
\section{Calculation of the rational solution, general scheme }
 1. We consider the natural
representation of the symmetric group $S_{n}$ (see [1]). By $(i;j)$
we denote the permutation which transposes $i$ and $j$ and preserves
all the rest. The
 $n{\times}n$ matrix which corresponds to $(i;j)$ is denoted by
 \begin{equation}
 P(i,j)=[p_{k,\ell}(i,j)], \quad (i{\ne}j).\end{equation}
 The elements $p_{k,\ell}(i,j)$ are equal to zero except
 \begin{equation}p_{k,\ell}(i,j)=1,\quad (k=i,\ell=j);\quad p_{k,\ell}(i,j)=1,\quad
(k=j,\ell=i),\end{equation}
 \begin{equation}p_{k,k}(i,j)=1,\quad
 (k{\ne}i,k{\ne}j).\end{equation} Now we introduce the matrices
 \begin{equation}
 P_{k}=P(1,k+1),
 \quad 1{\leq}k{\leq}n-1.\end{equation}
 and the $1{\times}(n-1)$ vector $e=[1,1, ...,1]$ ,
 and
 $n{\times}n$ matrix \begin{equation}
 S=\left[\begin{array}{cc}
   2-n & e\\
   e^{\tau} & 0 \\
 \end{array}\right].
 \end{equation}
 Using relations  (2.1)-(2.4) we deduce that
 \begin{equation}
 T=\sum_{k=1}^{k=n-1}P_{k}=(n-2)I_{n}+S.\end{equation}
 The eigenvalues of $T$ are defined by the equalities
 \begin{equation} \lambda_{1}=n-1,\quad \lambda_{2}=n-2,\quad
 \lambda_{3}=-1.\end{equation} The corresponding eigenvectors have
 the forms
 \begin{equation} V_{1}=\mathrm{col}[1,1,...,1],\quad V_{2}=\mathrm{col}[0,a_{1},...,a_{n-1}],\quad
V_{3}=\mathrm{col}[(n-1),-1,...,-1].\end{equation} 2. First we
investigate equation (1.1) , (1.2) in a neighborhood
 of $z=\infty$. Changing the variable $z=1/\xi$ we obtain
 \begin{equation}\frac{dV}{d\xi}=-{\rho}B(\xi)V(\xi),\end{equation}
 where
 \begin{equation}
 V(\xi)=W(1/\xi)=\sum_{p=m}^{\infty}{\xi}^{p}G_{p},\quad
 G_{m}{\ne}0,\quad |\xi|<r,\end{equation}
 \begin{equation}B(\xi)=\frac{1}{{\xi}^{2}}A(1/\xi)=T{\xi}+\sum_{p=0}^{\infty}{\xi}^{p}T_{p},
\quad |\xi|<r.\end{equation}Here the matrix $T$ is defined by
relations (2.5) and (2.6). In view of (1.2) and (2.11)  the
equalities
\begin{equation}
T_{p}=\sum_{k=1}^{n-1}P_{k}z_{k}^{p+1},\quad p{\geq}0\end{equation}
are valid.  In case $\rho=1$ we have the system
\begin{equation}
[(q+1)I_{n}+T]G_{q+1}=-\sum_{j+\ell=q}T_{j}G_{\ell},\quad
j{\geq}0,\quad q+1{\geq}-n+1.\end{equation} In case $\rho=-1$ we
obtain that
\begin{equation}
[(q+1)I_{n}-T]G_{q+1}=\sum_{j+\ell=q}T_{j}G_{\ell},\quad
j{\geq}0,\quad q+1{\geq}-1.\end{equation} We need the following
assertion (see [5]).\\
 \textbf{Proposition 2.1.}
\emph{Let $\rho=\pm{1}$.Then the matrix function $W(z)$ can be
written in the form}
\begin{equation}W(z)=\sum_{k=1}^{n-1}\frac{L_{k}}{z-z_{k}}+Q(z),\end{equation}
\emph{where $L_{k}$ are $n{\times}n$ matrices,
$Q(z)$ is  $n{\times}n$ matrix polynomial.}\\
Using relations (2.10) and (2.15) we deduce that
\begin{equation} \sum_{k=1}^{n-1}z_{k}^{p-1}L_{k}=G_{p},\quad
1{\leq}p{\leq}s,\quad \rho=\pm{1}.\end{equation} 3. Now we can
describe the method of calculating the rational solution $W(z)$ of
system (1.1),(1.2).\\
Step 1. From relation (2.13) we find $G_{p} (-n+1{\leq}p{\leq}n-1)$
in case $\rho=1$.  From relation (2.14) we find $G_{p}
(-1{\leq}p{\leq}n-1)$ in case $\rho=-1$.\\
Step 2. We find $Q(z)$ with the help of the formulas
\begin{equation}Q(z)=\sum_{q=0}^{n-1}z^{q}G_{-q},\quad
\rho=1,\end{equation}
\begin{equation}Q(z)=\sum_{q=0}^{1}z^{q}G_{-q},\quad
\rho=-1.\end{equation} Step 3. In view of (2.16) the matrices
$L_{k}$ are defined by the equality
\begin{equation}
\mathrm{col}[L_{1},L_{2},...,L_{n-1}]=M^{-1}\mathrm{col}[G_{1},G_{2},...,G_{n-1}].\end{equation}
 The block matrix $M$ has the form
\begin{equation}
M=\left[\begin{array}{cccc}
  I_{n} & I_{n} & ... & I_{n} \\
  z_{1}I_{n} & z_{2}I_{n} & ... & z_{n-1}I_{n} \\
  ... & ... & ... & ... \\
  z_{1}^{n-2}I_{n} & z_{2}^{n-2}I_{n} & ... & z_{n-1}^{n-2}I_{n}   \\
\end{array}\right].\end{equation}Together with matrix $M$ we consider
the Vandermonde matrix \\
$$M_{1}=\left[\begin{array}{cccc}
 1 & 1 & ... & 1 \\
  z_{1} & z_{2} & ... & z_{n-1} \\
  ... & ... & ... & ... \\
  z_{1}^{n-2} & z_{2}^{n-2} & ... & z_{n-1}^{n-2}   \\
\end{array}\right]$$\\
The following relation
\begin{equation}\mathrm{det}M_{1}=\prod_{i>j}(z_{i}-z_{j}){\ne}0\end{equation}
is true for Vandermonde matrix $M_{1}$. The matrix
$N=M^{-1}=[N_{i,j}]$ is defined by the relations
\begin{equation}N_{i,j}=A_{j,i}I_{n} ,\end{equation} where the numbers $A_{j,i}$
are adjoint minors of entries $m_{j,i}$ of the matrix $M_{1}$.\\
4. Now we use the described method in the case when $\rho=-1.$
 The $n{\times}1$ vectors
\begin{equation}G_{-1}=G_{0}=...=G_{n-2}=0,\quad
G_{n-1}=\mathrm{col}[1,1,...,1] \end{equation} satisfy  relations
(2.14). The vectors $L_{k}$ can be chosen in the following form
\begin{equation}L_{k}={\alpha}_{k}V_{k}.\end{equation} We note, that the vectors $V_{k}$
 are defined by relation (2.8).
 Now using formula
(2.19) we find the numbers \begin{equation}
\alpha_{k}=A_{n-1,k}/\mathrm{det}M_{1}=(-1)^{n+k-1}/[\prod_{i>k}(z_{i}-z_{k})\prod_{k>j}(z_{k}-z_{j})]
\end{equation}Relations
(2.23)-(2.25)
imply the following assertion.\\
 \textbf{Proposition 2.2.} \emph{The
vector function
\begin{equation}Y_{1}(z)=\sum_{k=1}^{k=n-1}\frac{L_{k}}{z-z_{k}}
\end{equation}
 is
the rational solution of system $(1.1)$ when $\rho=-1.$}\\
 In order to construct the new rational solutions $Y_{k}(z),\quad (2{\leq}k{\leq}n-1)$ of system (1.1)
 we introduce the vectors
 \begin{equation}G_{-1}=G_{0}=...=G_{n-3}=0,\quad
 G_{n-2}=\mathrm{col}[0,a_{1},a_{2},...,a_{n-1}],\end{equation}where
\begin{equation}
a_{1}+a_{2}+...+a_{n-1}=0.\end{equation} It follows from (2.14) that
\begin{equation}[(n-1)I_{n}-T]G_{n-1}=T_{0}G_{n-2}.\end{equation}
The right hand sight of (2.29) has the form
\begin{equation}T_{0}G_{n-2}=\mathrm{col}[m,m_{1},m_{2},...,m_{n-1}],\end{equation}
where
\begin{equation}m=\sum_{i=1}^{i=n-1}a_{i}z_{i},\quad
m_{k}=a_{k}\sum_{i{\ne}k}z_{i}.\end{equation} We use the following
equality
\begin{equation}T_{0}G_{n-2}=[m/(n-1)]V_{3}+\mathrm{col}[0,b_{1},b_{2},...,b_{n-1}],
\end{equation}where the vector $V_{3}$ is defined by (2.8) and
\begin{equation}b_{k}=m_{k}+m/(n-1).\end{equation}
It is easy to see that
\begin{equation}\sum_{i=1}^{i=n-1}b_{i}=0.\end{equation}Using
relations (2.7), (2.8) and (2.29), (2.32) we have
\begin{equation}
G_{n-1}=m/[n(n-1)]V_{3}+
\mathrm{col}[0,b_{1},b_{2},...,b_{n-1}]\end{equation} In view of
(2.19), (2.27)  and (2.35) the relations
\begin{equation}L_{k}=\ell_{k}=A_{n-2,k}G_{n-2}+A_{n-1,k}G_{n-1},\quad
1{\leq}k{\leq}n-1 \end{equation} are true.\\
\textbf{Remark 2.1.} The following relations
\begin{equation}A_{n-2,k}=-A_{n-1,k}\sum_{i{\ne}k}z_{i},\quad
1{\leq}k{\leq}n-1 ,\end{equation}
\begin{equation}
A_{n-1,k}=\prod_{i>j,\quad i,j{\ne}k}(z_{i}-z_{j})\end{equation} are valid.\\
Formulas (2.27), (2.35)- (2.38)
 imply the following assertion.\\
\textbf{Proposition 2.3.} \emph{The vector functions
\begin{equation}Y_{j}(z)=\sum_{k=1}^{k=n-1}\frac{\ell_{k,j}}{z-z_{k}},\quad
(2{\leq}j{\leq}n-1)
\end{equation}
are the linear independent rational solutions of system $(1.1)$ when
$\rho=-1.$ Here}
\begin{equation}a_{n-1,j}=-1,\quad a_{k,j}=0 \quad if\quad j{\ne}k+1,\quad
a_{k,k+1}=1.\end{equation}
 In order to construct the  rational solution $Y_{n}(z)$ of system (1.1)
 we introduce the vector
 \begin{equation}G_{-1}=V_{3}.\end{equation}The vector $V_{3}$ can be
 represented in the form
 \begin{equation}V_{3}=\sum_{i=1}^{i=n-1}N_{i},\end{equation}where
 the
 entries $n_{i,j}$ of the vectors $N_{i}$ are defined by the relations
\begin{equation} n_{i,1}=n_{i,i+1}=1;\quad n_{i,j}=-2/(n-2), \quad if
\quad j{\ne}1 \quad and \quad if \quad j{\ne}i+1.\end{equation}
According to (2.14) the equality
\begin{equation}
TG_{0}=-T_{0}G_{-1} \end{equation} holds. Let us calculate the
right-hand site of (2.44):
\begin{equation}T_{0}G_{-1}=(\sum_{k=1}^{s}P_{k}z_{k})(\sum_{i=1}^{i=n-1}N_{i})\end{equation}
Using relations
\begin{equation}P_{k}N_{s}=P_{s}N_{k}, \end{equation} we have
\begin{equation}T_{0}G_{-1}=T\sum_{i=1}^{i=n-1}z_{i}N_{i}.\end{equation}
It follows from formulas (2.44) and (2.47) that
\begin{equation}G_{0}=-\sum_{i=1}^{i=n-1}z_{i}N_{i}.\end{equation}
To find $G_{1}$ we use the relation
\begin{equation}
(I_{n}-T)G_{1}=T_{0}G_{0}+T_{1}G_{-1}.\end{equation} Equalities
(2.45)-(2.47)  imply that the right-hand site of (2.48) has the form
\begin{equation}
T_{0}G_{0}+T_{1}G_{-1}=\sum_{k>\ell}(z_{k}-z_{\ell})^{2}P_{k}N_{\ell}.
\end{equation}
 The vector $P_{k,\ell}=P_{k}N_{\ell}$ can be represented as the linear combination of the vectors
$V_{2}$ and
$V_{3}$:\begin{equation}P_{k,\ell}=U_{0}+U_{k,\ell}\end{equation}where
\begin{equation}
U_{0}=-2/[(n-1)(n-2)]V_{3} ,\end{equation}
\begin{equation}U_{k,\ell}=V_{2}=\mathrm{col}[0,a_{1,k,\ell},a_{2,k,\ell},...,a_{n-1,k,\ell}].
\end{equation}Here $a_{p,k,\ell}$ are defined by the relations
\begin{equation}a_{p,k,\ell}=-\frac{2}{n-2}-\frac{2}{(n-2)(n-1)}=-\frac{2n}{(n-2)(n-1)}\quad
if\quad p{\ne}k,\ell .\end{equation}If $p=k$ or $p=\ell$ we have
\begin{equation}
a_{p,k,\ell}=1-\frac{2}{(n-2)(n-1)}=\frac{n(n-3)}{(n-2)(n-1)}.\end{equation}
From formulas (2.49)-(2.55)  we deduce that
\begin{equation}
G_{1}=\sum_{k>\ell}(z_{k}-z_{\ell})^{2}G_{k,\ell},\end{equation}
where
\begin{equation}G_{k,\ell}=\frac{1}{2}U_{0}+\frac{1}{3-n}U_{k,\ell}.
\end{equation}
Let us represent the vector $G_{k,\ell}$ as linear combination of
the vectors $N_{i}$ (see (2.44)):
\begin{equation}G_{k,\ell}=\sum_{s=1}^{n-1}\beta_{s,k,\ell}N_{s}.\end{equation}According
to (2.58) we have
\begin{equation}\sum_{s=1}^{n-1}\beta_{s,k,\ell}=0,\end{equation}
\begin{equation}-\frac{2}{n-2}\sum_{s{\ne}p}\beta_{s,k,\ell}+\beta_{p,k,\ell}=\frac{2n}{(n-1)(n-2)(n-3)},\quad p{\ne}k,\ell.
\end{equation}If $p=\ell$ or $p=k$ then the equality
\begin{equation}-\frac{2}{n-2}\sum_{s{\ne}p}\beta_{s,k,\ell}+\beta_{p,k,\ell}=-\frac{n}{(n-1)(n-2)}\end{equation}
is true. From relations (2.59)-(2.61) we deduce that
\begin{equation} \beta_{p,k,\ell}=\frac{2}{(n-1)(n-2)},\quad
p{\ne}k,\ell,\end{equation} and
\begin{equation}\beta_{k,k,\ell}=\beta_{\ell,k,\ell}=-\frac{1}{(n-1)}.\end{equation}
In view of (2.16) the equality
\begin{equation} G_{1}=\sum_{s=1}^{n-1}L_{s},\end{equation}where
\begin{equation}L_{s}=\gamma_{s}N_{s}.\end{equation}In order to find
$\gamma_{s}$ we represent $G_{1}$ in the form (see (2.56)-(2.58)):
\begin{equation}G_{1}=\sum_{k>\ell}(z_{k}-z_{\ell})^{2}\sum_{s=1}^{n-1}\beta_{s,k,\ell}N_{s}\sum_{s=1}^{n-1}\beta_{s,k,\ell}N_{s}\quad.
\end{equation}
Hence we have
\begin{equation}\gamma_{s}=\sum_{k>\ell}(z_{k}-z_{\ell})^{2}\beta_{s,k,\ell}\quad.
\end{equation}Relations (2.43),(2.48) and (2.65), (2.67) imply the
following assertion.\\
 \textbf{Proposition 2.4.} \emph{The
vector function
\begin{equation}Y_{n}(z)=\sum_{k=1}^{k=n-1}\frac{L_{k}}{z-z_{k}}+zG_{-1}+G_{0}
\end{equation}
 is
the rational solution of system $(1.1)$ when $\rho=-1.$}\\
\textbf{Proposition 2.5.} \emph{Let us consider the case $S_{n}$
when $n{\geq}3,\quad \rho=-1.$ In this case system $(1.1)$ has the
fundamental rational solution of the form
\begin{equation}Y(z)=\sum_{k=1}^{k=n}c_{k}Y_{k}(z),\end{equation}
where $c_{k}$ are arbitrary constants}.\\
\textbf{Remark 2.2.} It is easy to see that the constructed
solutions $Y_{k}(z),\quad(1{\leq}k{\leq}n)$ are linearly independent
and satisfy the conditions of Proposition 1.2.\\
\textbf{Remark 2.3.} The fundamental rational solution can be
represented in the $ n{\times}n$ matrix form
\begin{equation}
W(z)=[Y_{1}(z),Y_{2}(z),...,Y_{n}(z)] .\end{equation}
 \textbf{Remark 2.4.} The explicit
solutions for the cases $S_{3}$ and $S_{4}$ were constructed in
papers [5] and [6]. The case $\rho=\pm{2}$ for $S_{3}$ was considered by A.Tydnyuk [7]\\
 We repeat the remark from paper [5].\\
\textbf{Remark 2.5.} Let the matrices $P_{k}$ be symmetric. If
system (1.1),(1.2) has a  fundamental rational $n{\times}n$ matrix
solution $W(z)$, when $\rho=m$, then this system has the
 fundamental rational $n{\times}n$ matrix solution
 $Y(z)=[W^{-1}(z)]^{\tau},$ when $\rho=-m.$
 (By the symbol $Q^{\tau}$ we denote the transposed matrix $Q$.)\\
 \textbf{Example 2.1.}(Without rational fundamental solution.)\\
 Let us consider the system
 \begin{equation}\frac{dW}{dx}=A(z)W,\quad z{\in}C,\end{equation}
 where $A(z)$ and $W(z)$ are $3{\times}3$ matrices. We suppose that
 $A(z)$ has the form
 \begin{equation}A(z)=\frac{m_{1}P_{1}}{z-z_{1}}+\frac{m_{2}P_{2}}{z-z_{2}}.\end{equation}
 Here $m_{1}$ and $m_{2}$ are integers,$z_{1}{\ne}z_{2}$, the
 matrices $P_{1}$ and $P_{2}$ are defined by relations (2.1)-(2.4).
 We introduce the matrix
 \begin{equation}T=\left[
                   \begin{array}{ccc}
                     0 & m_{1} & m_{2} \\
                     m_{1} & m_{2} & 0 \\
                     m_{2} & 0 & m_{1} \\
                   \end{array}
                 \right]\end{equation}
 The matrix $T$ has the following eigenvalues
 \begin{equation}\lambda_{1}=m_{1}+m_{2},\quad
 \lambda_{2,3}=\pm\sqrt{m_{1}^{2}-m_{1}m_{2}+m_{1}^{2}}.\end{equation}
 According to results of paper [4] the following statement is true.\\
 \textbf{Proposition 2.6.} \emph{If $\lambda_{2}$ and $\lambda_{3}$
 are not integer , then system  $(2.71)$ has no rational fundamental
 solution.} \\
 \textbf{Corollary 2.1.} \emph{In case $m_{1}=1,\quad m_{2}{\ne}0,1$
 system $(2.71)$ has no rational fundamental solution.}\\
 \emph{Proof.} The inequalities
 \begin{equation}(m_{2}-1)^{2}<1-m_{2}+m_{2}^{2},\quad
 m_{2}>1,\end{equation}
\begin{equation}m_{2}^{2}<1-m_{2}+(m_{2}-1)^{2},\quad
 m_{2}<0\end{equation}imply that $\lambda_{2,3}$ are not integer. In view
 of Proposition 2.6 the corollary is true.\\
 \section {Consistent System}
 In section 1 we have considered only one equation of the
 Knizhnik-Zamolodchikov system.The corresponding system has the form
 \begin{equation}
\frac{dW}{dz_{i}}={\rho}A_{i}(z_{1},z_{2}, ... ,z_{n})W,\quad
i=1,2,...n,\end{equation} where $A_{i}(z_{1},z_{2}, ... ,z_{n})$ and
$W(z_{1},z_{2}, ... ,z_{n})$ are matrices of $n{\times}n$ and
$n{\times}1$ order respectively. We suppose that
 $A_{i}(z_{1},z_{2},...,z_{n})$
 has the form
\begin{equation}
A_{i}(z_{1},z_{2},...,z_{n}) =\sum_{
k{\ne}i}\frac{P_{i,k}}{z_{i}-z_{k}},
\end{equation} where $z_{k}{\ne}z_{\ell}$ if $k{\ne}\ell$.\\
\textbf{Remark 3.1.}  The following properties of
$P_{i,j}$\\
1. $P_{i,j}=P_{j,i} ,$\\
2. $[P_{i,j}+P_{j,k},P_{i,k}]=0, \quad i,j,k \quad distinkt,$\\
3. $[P_{i,j},P_{k,\ell}]=0,\quad i,j,k,\ell \quad distinkt$\\
are well known.\\
From properties 1.-3. we deduce directly the following assertion.\\
\textbf{Proposition 3.1.} \emph{System (3.1), (3.2) is consistent.}
\begin{center}\textbf{References} \end{center}
1. Burrow M., Representation Theory of Finite Groups, Academic
Press,
1965.\\
2. Chervov A., Talalaev D., Quantum Spectral Curves, Quantum
Integrable Systems and the Geometric Langlands Correspondence,
arXiv:hep-th/0604128, 2006.\\
3. Etingof P.I., Frenkel I.B., Kirillov A.A. (jr.), Lectures on
Representation Theory and Knizhnik-Zamolodchikov Equations, Amer.
Math. Society, 1998.\\
4. Sakhnovich L.A., Meromorphic Solutions of Linear Differential
Systems, Painleve Type Functions, Operators and Matrices, v.1,
87-111, 2007.\\
5. Sakhnovich L.A., Rational Solutions of KZ equation, existence and
construction, arXiv:math-ph/0609067, 2006.\\
6. Sakhnovich L.A., Rational Solution of KZ Equation, Case $S_{4}$.\\
arXiv:math.CA/0702404, 2006 \\
7. Tydnyuk A. Rational Solution of the
KZ Equation
(Example), \\ arXiv:math/0162153, 2006.\\

\end{document}